\documentclass[12pt]{article}
\usepackage{amsmath}
 \usepackage{amssymb}
 \usepackage{array}
 \usepackage{geometry}
 \usepackage{graphicx}
 \usepackage{enumerate}
\newtheorem {thm}{Theorem}[section]
\newtheorem{lem}{Lemma}[section]

\newtheorem {cor}{Corrollory}[section]
\newtheorem {ex}{Example}[section]
\newtheorem {de}{Definition}[section]
\numberwithin{equation}{section} 
\begin{document}
\begin{center}
{\Large \textbf{A New Family of Fractional Renewal Processes }}\\
\end{center}

{\begin{center}  Jung Hun Han \footnote{Email :
 jhan176@gmail.com}\end{center}}

\begin{abstract}
Fractional renewal processes as a generalization of Poisson process are already
in the literature.
In this paper, by introducing a new concept of generalized density function,
the authors construct new fractional renewal processes in the
$\alpha$-fractional space
and show that it is another interesting and useful generalization of
Poisson process.
\end{abstract}

Keywords: L\'{e}vy density function, fractional renewal process, fractional
caculus, Mittag-Leffler function,
\\

AMS Subject Classifications :26A33, 33E12, 60G22, 60k05;

\section{\bf  Introduction}\label{s1}
In \cite[p.16]{gm2}, the authors states ``By time change via the inverse stable
subordinator the standard Poisson process is transformed to the fractional
Poisson process."
\vskip 0.1cm
In physics, the area of Beck-Cohen superstatistics is growing fast as can be
seen in \cite{bc1,h3,mh1,t1}, which belongs to Bayesian world in
Statistics. 
In \cite{bc1} and their other papers, the authors treat them under the two
categories, one of which is associated with a non-density function to find
useful information and interpretation of some physical systems in nature. In
statistical aspects, the pathway idea is one of areas growing dramatically in
this direction as shown in \cite{h3,mh1}. In mathematics,
we have similar situations as the Dirac delta generalized function shows up in
many places of pure mathematics even though it violates the true definition of a
function.
\vskip 0.1cm

In the theory of renewal process, Poisson process plays an enormous role since
it is the only one that comes from Bernoulli trial by the limit process. Many
similar generalizations have been attempted. In \cite{gm1}, the authors provide
a fractional renewal process, called a fractional Poisson process, and it
appears as a good generalization of the original Poisson process and Erlang
process. 

In \cite{h1}, the author insists that with respect to $\alpha$ the fractional
world has to be categorized and classified under the condition that its
L\'{e}vy structure be preserved where $\alpha$ is fixed and lies between $0$ and
$1$.  
 
In this paper, motivated by and combination of these ideas, the authors define a
generalized random variable along with its generalized distribution and apply
it to the world of the $\alpha$-fractional space. Furthermore the authors
introduce a family of new fractional renewal processes.

To see the role of L\'{e}vy structure, the correspondence via 
J-transformation from the author's paper \cite{h1} provides the following
\begin{equation}
 e^{-x} \longleftrightarrow t^{\alpha-1}
E_{\alpha,\alpha}(-t^{\alpha}).
\end{equation}
In this paper, there is an analogous correspondence
\begin{equation}
 \mbox{ Poisson process } \longleftrightarrow ~~\alpha\mbox{-fractional Poisson
process. }
\end{equation}
  
\section{$\alpha$-Fractional Space}
Note that $\alpha$ varies in $(0,1)$ throughout this paper and $\alpha$ is
fixed always.
\begin{de}
Let f(x) have $H$-function representation and be convergent
\cite{mh2, mhs1}. Then f(x) is said to be a function with L\'{e}vy
structure if its $H$-function representation has the factor
$\frac{\Gamma(-\frac{s}{\alpha}+\frac{1}{\alpha})}{\alpha\Gamma(1-s)}$
in the integrand. In short,
$\frac{\Gamma(-\frac{s}{\alpha}+\frac{1}{\alpha})}{\alpha\Gamma(1-s)}$
will be called L\'{e}vy structure.
\end{de}
Note that the Laplace transform of
$\displaystyle \frac{1}{2 \pi i}\oint_{L}
\frac{\Gamma(\frac{1}{\alpha} - \frac{s}{\alpha}  )}{\alpha
\Gamma(1-s) }x^{-s}ds, ~~ 1>Re(s)>0,~\alpha>0$ is $e^{-s^{\alpha}}$, which is
shown in \cite{ma1}. L\'{e}vy structure is named from this relation.

In \cite{ma1}, the author shows a process to lift a gamma density to a
generalized Mittag-Leffler density function by using statistical
techniques and it is in Example \ref{hhh31}.
\begin{ex}\label{hhh31}
Let $x$ be a simple exponential random variable with the density
function $f(x)=e^{-x}$. We attach the L\'{e}vy structure to $
E[x^{\frac{1}{\alpha}}]^{s-1} =\int^{\infty}_{0}
(x^{\frac{1}{\alpha}})^{s-1}e^{-x}dx$.
  Then
\begin{equation}
   x^{\alpha-1}E_{(\alpha,\alpha)} (- x^{\alpha})=\frac{1}{2\pi i} \oint_{L}
\frac{E[x^{\frac{1}{\alpha}}]^{s-1}\Gamma(-\frac{s}{\alpha}+\frac{1}{\alpha})x^{
-s}}
{\alpha\Gamma(1-s)}ds
  \end{equation}
by the residue theorem.
\end{ex}
In \cite{h1}, the author shows using Mellin transformation property how to lift
a function in the ordinary space to a corresponding function in the
$\alpha$-fractional(or $\alpha$-level) space in an analytic
way including statistical ones.

 Let $f(x)$ be a function, which does not have Levy structure and lives in the
ordinary space and $h(x)$ be the L\'{e}vy density function as a
kernel.
\begin{de}
 Define the J-transform of $f(x)$ by
\begin{eqnarray}
 J(f)(x) = \lim_{\gamma\rightarrow\infty}\int^{x}_{0}
\left(\frac{t}{x}\right)f_{2}\left(\left(\frac{x}{t}\right)^{\alpha}\right)
 \frac{t^{\alpha\gamma-1}\gamma^{\gamma}E^{\gamma}_{(\alpha,\alpha\gamma)}
(-\gamma t^{\alpha})}{t}dt
   \end{eqnarray}
   where ~$\alpha$ fixed in $0 < \alpha < 1$,$x\geq 0$ $f_{2}(x)=xf(x)$ and
$f(x)$ is integrable and continuous on
   the interval.
\end{de}

\begin{ex}\label{hhh16}
Let $ f(x)=e^{-x}$. Then
\begin{eqnarray*}
  J(f)(x) &=&x^{\alpha-1}E_{(\alpha,\alpha)} (- x^{\alpha}).
\end{eqnarray*}
\end{ex}
{\bf\large Ordinary space and $\alpha$-fractional space }\\
The concept or term of ordinary space looks strange. But we want the
ordinary space to be
generated by $\displaystyle \left\{0,1,x, x^{2},\ldots\right\}$,
namely the exponents of the variables
being non-negative integers only.  Then we obtain the
corresponding $\alpha$-fractional space 
generated by the elements spawned as J-transforms of
$\displaystyle \left\{0,1,x, x^{2},\ldots\right\}$. The basis for
$\alpha$-fractional space seems inherited from that of the ordinary space.
But from the piont of view of Probability theory, the following
correspondence looks natural:
\begin{equation}
 \left\{0,1,x,\frac{ x^{2}}{2!},\ldots\right\}
\longleftrightarrow
\left\{0,\frac{x^{\alpha-1}}{\Gamma(\alpha)},\frac{x^{2\alpha-1}}{
 \Gamma(2\alpha) } ,\frac{x^{3\alpha-1}}{\Gamma(3\alpha)},\ldots\right\}.
\end{equation}
Note that the factorial coefficients in front
of variables must appear as statistical quantities since Mellin transformation
and its family preserve and carry statistical measures.
Therefore every function will have the L\'{e}vy
structure in their $H$-function representations in the sense of the author
\cite{h1}. In the $\alpha$-fractional world, we can change to the alternative
definition of the Riemann-Liouville fractional derivative without the ordinary
derivative as follows:
\begin{equation}
 (D_{0}^{\alpha}\phi_{\alpha})(t)=\int_{0}^{t}\frac{(t-x)^{-\alpha-1}}{
\Gamma(-\alpha)} \phi_{\alpha}(x)dx.
\end{equation}
We emphasise that the main function $\phi_{\alpha}(t)$ is the solution of the
Reimann-Liouville fractional integral equation
$\phi_{\alpha}(t)-\frac{t^{\alpha-1}}{\Gamma(\alpha)}=
-(D^{-\alpha}_{0}\phi_{\alpha})(t)$ whereas $g_{1}(t)=e^{-t}$ is the solution of
$g_{1}(t)-1= -\int_{0}^{t}g_{1}(x)dx$.

\section{Standard Poisson Process and CTRW}
In \cite{gm1}, the authors describe the underlying theory and setting of renewal
processes (Poisson process is a special case of renewal process) and
continuous time random walk (CTRW). \\
{\bf\large Standard Poisson Process and Erlang Process}\\
Poisson process is suitable for some situations like modelling of counting the
number of arrivals of customers in a particular place in a given amount of time
or some physical system showing random behaviors such as a jump to the next
position. The probability for $n$ arrivals in the given time interval can be
calculated with $p(n,t_{1})$ for fixed $t_{1}$.
It is assumed that $\lambda=1$. Poisson process can be defined as
an infinite sequence $0=t_{0}<t_{1} <t_{2}
\cdots$ of events separated by i.i.d. (independent and identically distributed)
random waiting times $T_{j}=t_{j}-t_{j-1}$ of exponential distributions.
The defining characteristics of the Poisson process are  time
homogeneity, independence for defferent waiting time random variables and 
 infinitesimal interval probabilities.
These assumptions make $N(t)$ be a Poisson process with the
Poisson distribution with the parameter $t$ and the expectation
of $N(t)$ (the average number of
arrivals in the given time interval)
\begin{equation}
 p(n,t)=\frac{t^{n}}{n!}e^{-t}, ~~
m(t)=<N(t)>=t. 
\end{equation}
$m(t)$ is known as a renewal function.\\
 Its inverse process, which is called Erlang process, is designed to calculate
the probability for the time when the $n$-th arrival has just happened with its
Erlang density function
$\displaystyle q_{n}(t)=\frac{t^{n-1}}{(n-1)!}e^{-t}, t\geq 0$ and its
cumulative distribution $Q_{n}(t)$.\\
In the probailistic language, the following property is called memoryless:
\begin{equation}
 \mathit{P}\{X>t+x\}= \mathit{P}\{X>t\}\mathit{P}\{X>x\}.
\end{equation}
Since the exponential function satisfies this property only, Poisson process is
expected to have memoryless property. But the Mittag-Leffler density function,
which is the main function in this paper, does not have this property but
  the property of long-term-memory type.  \\
{\bf\large Continuous Time Random Walk}\\
In the theory of CTRW, by using Dirac delta functions, one is able to
develop random walk model in discrete space as a special case of CTRW which is
known as a compound Poisson process associated with random walk. It can be
considered that time variable and space variable vary in the positive real line:
 choose $w(x)=\delta(x-1), \phi_{\alpha}(t)=t^{\alpha-1}
E_{\alpha,\alpha}(-t^{\alpha}),x\geq 0, t> 0$. Notice that as indicated in
\cite{gm1}, we can ignore the possible delta peak at the origin of the time
line by taking $t>0$. So the Cox-Weiss series for CTRW can be used and the
Laplace-Laplace solution shall be taken to be:
\begin{equation}\label{r8}
 \widetilde{\widetilde{p}}_{\alpha}(k,s):= \widetilde
{\Psi}_{\alpha}(s)\sum_{n=0}^{\infty} (
\widetilde{\phi}_{\alpha}(s) \widetilde{w}(k))^{n} = \widetilde
{\Psi}_{\alpha}(s) \frac{ 1}{ 1-
\widetilde{\phi}_{\alpha}(s) \widetilde{w}(k)},~
p_{\alpha}(x,t):=\sum_{n=0}^{\infty} P_{n,\alpha}(t)
\delta(x-n)
\end{equation}
where $\Psi_{\alpha}(t)$ is the marginal distribution
of the density
function $\phi_{\alpha}(t)$, which is defined in section \ref{5}.\\

In section \ref{5}, we shall find out some analogues of the integral equation
of the CTRW and the Kolmogorov-Feller equation (the master equation of the
compound Poisson process).

\section{Abel-Volterra Equation of the Second Kind}\label{4}
\begin{equation}
  \label{hhh1} f(x)-g(x)= - c(D^{-\alpha}_{0}f)(x)
  \end{equation}
where $g(x)$ is a general integrable function on the finite
interval $[0,b]$ and $0<\alpha<1$. For detailed history, theory and
applications, see \cite{gv1}. This integral equation is known as
Abel-Volterra integral equation of the second kind.
The solution of the Abel-Volterra integral equation is
\begin{eqnarray}\label{hhh10}
f(x)&=& g(x)-c \int^{x}_{0}(x-t)^{\alpha-1}E_{\alpha,\alpha}(-c(x-t)^{\alpha}
)g(t)dt.
\end{eqnarray}
The following table from the author's paper \cite{h1} gives certain densities
and their L\'{e}vey structures.\\
\textbf{\emph{functions without Levy structure}}\\
\begin{tabular}{|c|c|c|}
  \hline
  $g(x)$ &$\rightarrow$  & $f(x)$ \\
  \hline
  $1$ &  & $E_{\alpha}(-cx^{\alpha})$ \\
   \hline
   $x$&  &$xE_{\alpha,2}(-cx^{\alpha})$  \\
   \hline
  $e^{-x}$ &  &
$\sum^{\infty}_{k=0}(-x)^{k}E_{\alpha,k+1}(-cx^{\alpha})$ \\
   \hline
   $\frac{x^{\mu-1}}{\Gamma(\mu)}$ &  &
$x^{\mu-1}E_{\alpha,\mu}(-cx^{\alpha})$ \\
   \hline
   $x^{\mu-1}E^{\gamma}_{(\alpha,\mu)}(-cx^{\alpha})$& 
&$x^{\mu-1}E^{\gamma+1}_{\alpha,\mu}(-cx^{\alpha})$  \\
   \hline
\end{tabular}
\\

\vskip 0.2cm
\textbf{\emph{functions with Levy structure}}\\
\begin{tabular}{|c|c|c|}
  \hline
  $g(x)$ &$\rightarrow$  & $f(x)$ \\
  \hline
   $\frac{x^{\alpha-1}}{\Gamma(\alpha)}$ &  &
$x^{\alpha-1}E_{\alpha,\alpha}(-cx^{\alpha})$ \\
   \hline
   $x^{\alpha-1}E^{\gamma}_{\alpha,\alpha}(-cx^{\alpha})$& 
&$x^{\alpha-1}
   E^{\gamma+1}_{\alpha,\alpha}(-cx^{\alpha})$  \\
   \hline
  \end{tabular}
  \vskip 0.1cm

So when $c= (1-\widetilde{w(k)})$ and $g(x)=1$, the
Abel-Volterra integral equation becomes 
\begin{equation}
 f(t)-1=-(1-\widetilde{w(k)})(D_{0}^{-\alpha} f)(t).
\end{equation}
Then $E_{\alpha}(-(1-\widetilde{w(k)})t^{\alpha})$, which appears in \cite{gm1},
becomes the solution of the
above  Abel-Volterra integral equation.
As suggested in \cite{h1}, we consider an Abel-Volterra integral
equation with Levy structure preserved as follows:
Set $g(t)= \frac{t^{\alpha-1}}{\Gamma(\alpha)}$ and
$N_{0}=1$ where $c =(1-\widetilde{w(k)})$. 
Then the Abel-Volterra integral equation becomes 
\begin{equation}
 f(t)-\frac{t^{\alpha-1}}{\Gamma(\alpha)}= -
c(D^{-\alpha}_{0}f)(t).
\end{equation}
We have the following solution
\begin{equation}\label{430}
 f(t)=t^{\alpha-1}
   E_{\alpha,\alpha}(-ct^{\alpha}).
\end{equation}
Define $\widetilde{p}_{\alpha}(k,t):= t^{\alpha-1}E_{\alpha,\alpha}
(-(1-w(k))t^{\alpha})$.
From the Laplace-pair table, we have
\begin{eqnarray}
\widetilde{\widetilde{p}}_{\alpha}(k,s)&=& \frac{1}{s^{\alpha}}
\left(1+\frac{(1-\widetilde{w(k)})}{s^{\alpha}}\right)^{-1} =
\frac{1}{(1+s^{\alpha})}
\frac{1}{\left(1-\frac{\widetilde{w(k)}}{1+s^{\alpha}}\right)}
\\
&=&\sum_{n=0}^{\infty} 
\widetilde{\phi}^{1+n}_{\alpha}(s)\widetilde{w(k)}^{n} \label{a6}
\end{eqnarray}
When $\widetilde{w(k)}=e^{-k}$,  we have
\begin{equation}
 p_{\alpha}(x,t) =\sum_{n=0}^{\infty}
\left(  \phi^{*(1+n)}_{\alpha}\right)(t)
\delta(x-n).
\end{equation}

\section{A Fractional Renewal Process of Order $\alpha$}\label{5}

\begin{de}
By {\bf an asymptotic limit}, we mean that if a density is
asymptotically equivalent to a power function of the form 
$\frac{x^{\alpha-1}}{\Gamma(\alpha)}$ for a fixed $\alpha$ with $|\alpha|<1$
when its parameter or variable approachs to infinity, then the power function is
said to be the {\bf asymptotic limit} of the density function. 
\end{de}
\begin{ex}
  The Gaussian normal density has the asymptotic limit as $t\rightarrow\infty$  
\begin{equation}
\frac{t^{-\frac{1}{2}}}{\Gamma(\frac{1}{2})}
\exp\left(-\frac{x^{2}}{t}\right)\longrightarrow
\frac{t^{1/2-1}}{\Gamma\left(\frac{1}{2} \right)}.
\end{equation}
\end{ex}
\begin{ex}
The L\'{e}vy density function of order $\frac{1}{2}$ has the asymptotic
limt
\begin{equation}
 \frac{at^{-\frac{3}{2}}}{2\Gamma(\frac{1}{2})}
\exp\left(-\frac{a^{2}}{t}\right)=\frac{at^{-\frac{3}{2}}}{-\Gamma(-\frac{1}{2})
}
\exp\left(-\frac{a^{2}}{t}\right)\longrightarrow
\frac{t^{-\frac{1}{2} -1}}{\Gamma(-\frac{1}{2})} 
\end{equation}
 as $a$ goes to $-1$ and $t$ tends to infinity.
\end{ex}
Note that for the Green functions appearing in the above examples, we refer to
 \cite{mpg1}, also refer to \cite{f1} for the systematic treatment of L\'{e}vy
density functions and the author provides the explcit formulae of L\'{e}vy
stable density functions of the type $\displaystyle
\left(\frac{p}{q}\right)^{\frac{l_{2}}{l_{1}}}$ in \cite{h2}.
\begin{de}
  An asymptotic limit is said to be {\bf a generalized density function}
if the limit is from a density function and it makes sence when the function is
inside an integral. In fact, it is not a density function in general and it may
be called a supplimentary density function following the spiritual concept of
Dirac delta generalized function. Accordingly {\bf a generalized random
variable} is to be
defined with its generalized density function.
\end{de}
{\bf Remark.} We have generalized random variables of order $\alpha$ with
$|\alpha|<1$. When $\alpha$ is equal to zero, then the power function become the
Dirac delta function in the sense of Gel'fond and Shilov \cite{gs1}. Here we
assume that we are able to obtain asymptotic limits of all of $\alpha$. The
generalized density functions of the form $\frac{x^{\alpha-1}}{\Gamma(\alpha)},
0<\alpha<1$, plays a role as the Riemann-Liouville fractional integral and the
generalized density functions of the form
$\frac{x^{-\alpha-1}}{\Gamma(-\alpha)},
0<\alpha<1$, plays a role as the Riemann-Liouville fractional derivative in the
$\alpha$-fractional world.
\vskip 0.2cm

Let $X_{1}$ be the generalized random variable of order $\alpha$ fixed with
$0<\alpha<1$ and $S_{1}$ be the Mittag-Leffler random variable of the same order
$\alpha$. Then $X_{1}$ has the corresponding generalized density function of the
form $\displaystyle \frac{t^{\alpha-1}}{\Gamma(\alpha)}$ and $S_{1}$ does the
density function of the form $\phi_{\alpha}(t)=t^{\alpha-1}
E_{\alpha,\alpha}(-t^{\alpha})$.
We construct a world $\Omega$ of random variables with marginal
distributions obtained by considering $X_{1}+S_{1}$ from the density function.
Then we have the following relations:
\begin{eqnarray}
\Phi_{\alpha}(t)&=&(~_{0}D_{t}^{-\alpha}\phi_{\alpha})(t)
=\int_{0}^{t}\frac{(t-t_{1})^{\alpha-1}}{\Gamma(\alpha)}\phi_{\alpha}(t_{1})dt_{
1}, \\
\Psi_{\alpha}(t)&=& (~_{t}W_{\infty}^{\alpha}\phi_{\alpha})(t)=\int_{t}^{\infty}
\frac{(t_{1}-t)^{\alpha-1}}{\Gamma(\alpha)}\phi_{\alpha}(t_{1})dt_{1}
= \frac{t^{\alpha-1}}{\Gamma(\alpha)}-\Phi_{\alpha}(t).\label{k100}
\end{eqnarray}
It can be called $\Phi_{\alpha}(t):={\mathit P}(T\leq t)$ the failure
marginal probability and $\Psi_{\alpha}(t):={\mathit P}(T>t)$ the survival
marginal Probability. Note that the subscript $\alpha$ is used to indicate that
the
function lives in the $\alpha$-fractional space.

\begin{lem}\label{l1}
The following equalities are true:
\begin{equation}
  \phi_{\alpha}(t)=-(~_{0}D_{t}^{\alpha}\Psi_{\alpha})(t),
~~~~~~~~~\Psi_{\alpha}(t)=\phi_{\alpha}(t).
\end{equation}

\end{lem}
{\bf Proof}
Use the alternative definition for $D_{0}^{\alpha}$.$\square$\\
This lemma enables us to define a generalization of the Poisson process in the
$\alpha$-fractional space. Note that the process
$~_{0}D_{t}^{\alpha}\Psi_{\alpha}(t)$ is the reverse process from the marginal
distribution to the density function.

In this section, we frequently calculate $\phi^{*n}_{\alpha}$ for $n\in
\mathbb{N}$. So we
present the following formula:
it is known that $\displaystyle \frac{(n)_{k}}{k!},~ k=0,1,2,3,\ldots,$
generates the sequence of
counting numbers when $n=2$, that of triangular numbers when $n=3$, that of
tetrahedral numbers when $n=4$ and so on, which all come from the diagonal
sequences of Pascal's triangle. 
The following equalities are necessary:
\begin{eqnarray}
\sum_{k=0}^{\infty}t^{k}\sum_{n=0}^{\infty}t^{n}&=& \sum_{k=0}^{\infty}
(\sum_{n=0}^{k}1)t^{k}= \sum_{k=0}^{\infty}(k+1)t^{k} \label{426}
\\
\sum_{k=0}^{\infty}(k+1)t^{k}\sum_{n=0}^{\infty}t^{n} &=&
\sum_{k=0}^{\infty}\frac{(k+1)(k+2)}{2} t^{k} \label{427}
\\
\sum_{k=0}^{\infty}\frac{(k+1)(k+2)}{2} t^{k} \sum_{n=0}^{\infty}t^{n}
&=& \sum_{k=0}^{\infty}\frac{(k+1)(k+2)(k+3)}{3!} t^{k} \label{428}
\\
&&\cdots.
\end{eqnarray}
The coefficients of the above series (\ref{426}), (\ref{427}), (\ref{428}) and
so on also generate the diagonal sequences of Pascal's triangle. It can be
verified using the
properties of the diagonal sequences of Pascal's triangle. Therefore the
sequences generated by
$\displaystyle \frac{(n)_{k}}{k!},~ k=0,1,2,3,\ldots$ are equal to those
obtained from the
coefficients of the above series with the correspondence that $\displaystyle
n=2\longleftrightarrow (\ref{426}),~ n=3\longleftrightarrow (\ref{427}),~
n=4\longleftrightarrow (\ref{428})$ and so on . By using these properties, we
obtain
\begin{equation}
 \phi^{*n}_{\alpha}=t^{n\alpha -1}E^{n}_{\alpha,n\alpha}(-t^{\alpha}) =
\frac{t^{n\alpha
-1}}{(n-1)!} \sum_{k=0}^{\infty}\frac{\Gamma(n+k)}{\Gamma( \alpha n +\alpha k)}
\frac{(-1)^{k}t^{\alpha k}}{k!}.\label{413}
\end{equation}

\subsection{Mittag-Leffler function $t^{\alpha-1}
   E_{\alpha,\alpha}(-(1-\widetilde{w}(k))t^{\alpha}) $}
{\large\bf Cox-Weiss Series and Analogue of Montroll-Weiss Equation}\\
We have the following Cox-Weiss series
\begin{equation}
\widetilde{\widetilde{p}}_{\alpha}(k,s)=
\widetilde{\Psi_{\alpha}}(s)\sum_{n=0}^{\infty}
\widetilde{\phi_{\alpha}(s)}^{n}\widetilde{w(k)}^{n}.
\end{equation}
Since $\displaystyle \widetilde{\Psi_{\alpha}}(s) = \frac{1-
\widetilde{\phi_{\alpha}}(s)}{s^{\alpha}}$, we obtain the
analogue of the famous Montroll-Weiss equation from the Cox-Weiss series. 
\begin{equation}
\widetilde{\widetilde{p}}_{\alpha}(k,s)= \frac{1-
\widetilde{\phi_{\alpha}}(s)}{s^{\alpha}}\sum_{n=0}^{\infty}
\widetilde{\phi_{\alpha}(s)}^{n}\widetilde{w(k)}^{n}.
\end{equation}
The Laplace transform of $\phi_{\alpha}(s)$ is $\frac{1}{1+s^{\alpha}}$. Hence
we get
\begin{equation}
\widetilde{\widetilde{p}}_{\alpha}(k,s)=
\widetilde{\phi_{\alpha}(s)}\sum_{n=0}^{\infty}
\widetilde{\phi_{\alpha}(s)}^{n}\widetilde{w(k)}^{n} \label{f4}.
\end{equation}
which is equal to (\ref{a6}).\\
{\bf\large  Renewal Function}\\
Following the procedure in \cite{gm1}, when
$\widetilde{w}(k)=e^{-k}$, we get 
\begin{equation}
 p_{\alpha}(x,t)=\sum_{n=0}^{\infty} P_{n,\alpha}(t)\delta(x-n) =
\sum_{n=0}^{\infty} (\Psi_{\alpha} *
\phi^{*n}_{\alpha})(t)\delta(x-n)
\end{equation}
and 
\begin{equation}
 m_{\alpha}(t)=-\frac{\partial}{\partial k}\tilde{p}_{\alpha}(k,t)|_{k=0}
=\left(
\sum_{n=0}^{\infty} nP_{n,\alpha}(t) e^{-nk} \right)|_{k=0}= \sum_{n=0}^{\infty}
nP_{n,\alpha}(t).
\end{equation}
By applying this and $\displaystyle \sum_{n=0}^{\infty}nz^{n}=
\frac{z}{(1-z)^{2}},|z|<1$, we obtain
\begin{equation}
 \tilde{m}_{\alpha}(s)=\widetilde{\Psi_{\alpha}(s)}\sum_{n=0}^{\infty} n
\widetilde{\phi_{\alpha}(s)}^{n}
=\frac{1-\widetilde{\phi_{\alpha}(s)}}{s^{\alpha}} 
\left(\frac{\widetilde{\phi_{\alpha}(s)}}{(1-\widetilde{\phi_{\alpha}
(s)})^{2}}\right)=
\frac{\widetilde{\phi_{\alpha}(s)}}{s^{\alpha} (1-\widetilde{\phi_{\alpha}
(s)})}.
\end{equation}
Therefore we have derived the reciprocal pair of the
relationships in the Laplace domain
\begin{equation}
 \widetilde{m_{\alpha}(s)}=\frac{\widetilde{\phi_{\alpha}(s)}}{s^{\alpha}
(1-\widetilde{\phi_{\alpha} (s)})}
, ~~~~ \widetilde{\phi_{\alpha}(s)}=
\frac{s^{\alpha}\widetilde{m_{\alpha}(s)}}{1+
s^{\alpha}\widetilde{m_{\alpha}(s)}}.
\end{equation}
Hence the renewal equation is 
\begin{equation}
 m_{\alpha}(t)=(D_{0+}^{-\alpha}\phi_{\alpha})(t)
+(m_{\alpha}*\phi_{\alpha})(t)
\end{equation}
On the other hand, we have the expression
\begin{equation}\label{f15}
 \widetilde{p}_{\alpha}(k,t)=t^{\alpha-1}E_{(\alpha,\alpha)}(-(1-e^{-k})t^{
\alpha} ).
\end{equation}
From this, we derive the following
\begin{equation}
 m_{\alpha}(t)=-\frac{\partial}{\partial k}\widetilde{p}_{\alpha}(k,t)|_{k=0} =
\frac{t^{2\alpha-1}}{\Gamma(2\alpha)}, ~~0<\alpha<1.
\end{equation}
\vskip 1cm
 We have a well-known infinite system of differential-difference
equations for the Poisson process with intensity $\lambda >0$ and $t\geq 0$,
\begin{equation}
 P_{0}(t)= e^{-\lambda t}, \frac{d}{dt}P_{n}(t)
=\lambda(P_{n-1}(t)-P_{n}(t)), 
n\geq1,
\end{equation}
with initial conditions $P_{n}(0)=0, n=1,2,3,\ldots$, which can be
used to define the Poisson process. We give an analogous system of
fractional
differential-difference equations for the fractional poisson process in the
$\alpha$-fractional  space.\\
In \cite{p1}, the following Laplace pair is provided:
\begin{equation}
 \frac{t^{\alpha(n+1)-1}}{n!}
E_{(\alpha,\alpha)}^{(n)}(-t^{\alpha}) \longleftrightarrow
\frac{1}{(1+s^{\alpha})^{n+1}}.
\end{equation}
So we get 
\begin{equation}
 \widetilde{p}_{\alpha}(k,t)=\sum_{n=0}^{\infty} \frac{t^{\alpha -1} t^{\alpha
n}}{ n!}
E_{(\alpha,\alpha)}^{(n)}(-t^{\alpha}) \widetilde{w}^{n}(k).
\end{equation}

{\large\bf Counting Probabilities}\\
Therefore the counting probabilities are
\begin{equation}
 P_{n,\alpha}(t)=\mathit{P}\{N(t)=n\}= \frac{t^{\alpha -1} t^{\alpha n}}{ n!}
 E_{\alpha,\alpha}^{(n)}(-t^{\alpha}).
\end{equation}
Or directly
\begin{equation}
 P_{n,\alpha}(t)=\phi^{*(n+1)}_{\alpha}(t)= t^{(n+1)\alpha
-1}E^{n+1}_{\alpha,(n+1)\alpha}(-t^{\alpha}) = \frac{t^{(n+1)\alpha
-1}}{n!} \sum_{k=0}^{\infty}\frac{\Gamma(n+1+k)}{\Gamma( \alpha n+\alpha +\alpha
k)}
\frac{(-1)^{k}t^{\alpha k}}{k!}
\end{equation}
{\bf\large Fractional Differential-Difference Equations}\\
Then we have the following analogous differential-difference Equations
\begin{eqnarray}
\widetilde{P}_{n,\alpha}(s)&=& \frac{1}{(1+s^{\alpha})^{n+1}}
\\
(1+s^{\alpha}) \widetilde{P}_{n,\alpha}(s)&=&\widetilde{P}_{n-1,\alpha}(s)
\\
s^{\alpha}\widetilde{P}_{n,\alpha}(s)&=& \widetilde{P}_{n-1,\alpha}(s)-
\widetilde{P}_{n,\alpha}(s)
\\
P_{0,\alpha}(t)&=&t^{\alpha-1}E_{\alpha,\alpha}(-t^{\alpha}),
(D_{0+}^{\alpha} P_{n,\alpha})(t)= P_{n-1,\alpha}(t)-P_{n,\alpha}(t).
\end{eqnarray}
where $D_{0+}^{\alpha}$ is defined in the appendix.\\
\noindent {\large\bf Erlang Densities}\\
Since $q_{n,\alpha}(t)= \phi^{*n}_{\alpha}(t)$, it can be calculated directly. 
From the formula (\ref{413}), it can be easily checked that
\begin{equation}
 q_{n,\alpha}(t)=t^{n\alpha -1}E^{n}_{\alpha,n\alpha}(-t^{\alpha}) =
\frac{t^{n\alpha
-1}}{(n-1)!} \sum_{k=0}^{\infty}\frac{\Gamma(n+k)}{\Gamma( \alpha n +\alpha k)}
\frac{(-1)^{k}t^{\alpha k}}{k!}.
\end{equation}
{\large\bf Fractional Integral Equation of the Cotinuous Time Random Walk}\\
We can derive the fractional integral equation of the CTRW as follows:
\begin{eqnarray*}
 \widetilde{\widetilde{p}}_{\alpha}(k,s)&=&\frac{\widetilde{\phi}_{\alpha}(s)}{
1-\widetilde{\phi}_{\alpha}(s)\widetilde{w}(k)}
\\
 \widetilde{\widetilde{p}}_{\alpha}(k,s)-\widetilde{\phi}_{\alpha}(s)\widetilde{
w}(k)
\widetilde{\widetilde{p}}_{\alpha}(k,s)&=& \widetilde{\phi}_{\alpha}(s)
\\
\widetilde{\widetilde{p}}_{\alpha}(k,s) &=&
\widetilde{\phi}_{\alpha}(s)+\widetilde{\phi}_{\alpha}(s)\widetilde{w}(k)
\widetilde{\widetilde{p}}_{\alpha}(k,s)
\\
p_{\alpha}(x,t)&=&\Psi_{\alpha}(t)\delta(x) + \int_{0}^{t}
\phi_{\alpha}(t-t_{1})dt_{1} \int_{0}^{x}
w(x-x_{1}) p_{\alpha}(x_{1},t_{1})dx_{1}.
\end{eqnarray*}
{\large\bf Fractional Version of Kolmogorov-Feller Equation}
\begin{eqnarray}
\widetilde{\widetilde{p}}_{\alpha}(k,s)&=&\frac{\widetilde{\phi}_{\alpha}(s)}{
1-\widetilde{\phi}_{\alpha}(s)\widetilde{w}(k)}
\\
\widetilde{\widetilde{p}}_{\alpha}(k,s)&=& \frac{1}{1+s^{\alpha}} \frac{1}{
1-\widetilde{\phi}_{\alpha}(s)\widetilde{w}(k)}
\\
(1+s^{\alpha})\widetilde{\widetilde{p}}_{\alpha}(k,s)&=& 
\frac{1}{1-\widetilde{\phi}_{\alpha}(s)\widetilde{w}(k)}
\\
s^{\alpha}\widetilde{\widetilde{p}}_{\alpha}(k,s)&=&-\widetilde{\widetilde{p
}}_{\alpha}(k,s)+ 
\frac{1}{1-\widetilde{\phi}_{\alpha}(s)\widetilde{w}(k)}-1+1
\\
s^{\alpha}\widetilde{\widetilde{p}}_{\alpha}(k,s)&=&-\widetilde{\widetilde{p
}}_{\alpha}(k,s)+ \widetilde{w}(k)\widetilde{\widetilde{p
}}_{\alpha}(k,s)+1
\\
(D_{0+}^{\alpha} p_{\alpha}(x,t))(t)&=& -p_{\alpha}(x,t)
+\int_{0}^{x}w(x-x_{1})p_{\alpha}(x_{1},t)dx_{1} + \delta(x)\delta(t)
\\
(D_{0+}^{\alpha} p_{\alpha}(x,t))(t)&=& -p_{\alpha}(x,t)
+\int_{0}^{x}w(x-x_{1})p_{\alpha}(x_{1},t)dx_{1}
\\
(D_{0+}^{\alpha} p_{\alpha}(x,t))(t)&=&  -p_{\alpha}(x,t)
+\int_{0}^{x}\delta(x-x_{1})p_{\alpha}(x_{1},t)dx_{1}.
\end{eqnarray}
Notice that we have taken $t>0$ to avoid the problem at $t=0$. Therefore
$\delta(t)=0$ which explains the above equation.

\subsection{Another Intereting Density of Mittag-Leffler Type}
Take a look at $\phi^{*l}_{\alpha}(t)=t^{l\alpha
-1}E^{l}_{\alpha,l\alpha}(-t^{\alpha})$. The $H$-function representation of
this function is\\ $\displaystyle \frac{1}{2 \pi i}\oint_{L}
\frac{\Gamma(l-\frac{1}{\alpha} + \frac{s}{\alpha})}{\Gamma(l)  }
\frac{\Gamma(\frac{1}{\alpha} -\frac{s}{\alpha})}{\alpha\Gamma(1-s)} t^{-s}ds$
which has the L\'{e}vy structure. This function appears
in the literature and is called a generalized Mittag-Leffler density, see
\cite[p.28]{hms2}. And its Laplace transform is $\displaystyle
\frac{1}{(1+s^{\alpha})^{l}}$. The most interesting part is its relationship
with the fractional Poisson process which we have developed in section 5.1.

Define $\widetilde{p}_{\alpha,l}(k,t):=t^{l\alpha
-1}E^{l}_{\alpha,l\alpha}(-(1-\widetilde{w(k)})t^{\alpha})$. Then we get
\begin{eqnarray}
\widetilde{\widetilde{p}}_{\alpha,l}(k,s)&:=&\widetilde
{ \widetilde{p}}_{ \alpha}^{l}(k,s)=
\underbrace
{\left(\frac {
\widetilde { \phi_ {\alpha} (s)}}{1-\widetilde{\phi_{\alpha}
(s)}\widetilde{w(k)}}\right)\cdots \left(\frac{
\widetilde { \phi_ {\alpha} (s)}}{1-\widetilde{\phi_{\alpha}
(s)}\widetilde{w(k)}}\right) }_{l \mbox{ times}}
\\
&=&\widetilde{\phi_{\alpha}
(s)}^{ l}
 \sum_{n=0}^{\infty}\frac{(l)_{n}}{n!}
\widetilde{\phi_{\alpha}(s)}^{n}\widetilde{w(k)}^{n}
=\widetilde{\phi_{\alpha}(s)
}^{ l}
 \sum_{n=0}^{\infty}{n+l-1 \choose n}
\widetilde{\phi_{\alpha}(s)}^{n}\widetilde{w(k)}^{n}. 
\end{eqnarray}
The corresponding analogous Cox-Weiss series is :
\begin{equation}
 p_{\alpha,l}(x,t)=p_{\alpha}^{*l}(x,t)= \sum_{n=0}^{\infty}
\frac{(l)_{n}}{n!}(\Psi_{\alpha}^{*l} *\phi^{*n}_{\alpha})(t)w^{*n}(x).
\end{equation}
And we obtain 
\begin{equation}
 \widetilde{p}_{\alpha,l}(k,t)=\sum_{n=0}^{\infty
}\frac{ (l)_{n}}{n!}\left(\frac{t^{\alpha l -1} t^{\alpha n}}{(n+l-1)!}
 E_{(\alpha,\alpha)}^{(n+l-1)}(-t^{\alpha})\right)
\widetilde{w}^{n}(k).
\end{equation}
{\large\bf Counting Probabilities}\\
Therefore the counting probabilities are
\begin{equation}
 P_{n,\alpha,l}(t)=\mathit{P}\{N_{1}(t)+N_{2}(t)+\cdots+ N_{l}(t)
=n\}=\frac{ (l)_{n}}{n!}\frac{t^{\alpha l -1} t^{\alpha n}}{(n+l-1)!}
 E_{(\alpha,\alpha)}^{(n+l-1)}(-t^{\alpha})
\end{equation}
where $N_{j}(t)$ is I.I.D for $j=1\ldots l$. Or directly
\begin{eqnarray}
 P_{n,\alpha,l}(t)&=& \frac{(l)_{n}}{n!}
t^{(n+l)\alpha
-1}E^{n+l}_{\alpha,(n+l)\alpha}(-t^{\alpha})
\\
 &=&{n+l-1\choose n}\frac{1}{(n+l-1)!}
t^{(n+l)\alpha -1}
\sum_{k=0}^{\infty}\frac{\Gamma(n+l+k)}{\Gamma( \alpha
(n+l) +\alpha k)}
\frac{(-1)^{k}t^{\alpha k}}{k!}
\end{eqnarray}
{\large\bf Renewal Function}\\
For the renewal function and the density function, we obtain
\begin{equation}
 \widetilde{m}_{\alpha,l}(s)=\widetilde{\Psi_{\alpha}(s)} ^{l}
 \sum_{n=0}^{\infty} \frac{(l)_{n}}{n!}n\widetilde{\phi_{\alpha}(s)}^{n}
=\frac{l\widetilde{\phi_{\alpha}(s)}}{s^{\alpha l}(1-
\widetilde{\phi}_{\alpha}(s))},~~\widetilde{\phi_{\alpha}(s)}=
\frac{s^{\alpha l}\widetilde{m}_{\alpha,l}(s)}{(l+s^{\alpha l}
\widetilde{m}_{\alpha,l}(s))}
\end{equation}
since $\displaystyle \sum_{n=0}^{\infty} \frac{(l)_{n}nz^{n}}{n!}=
\frac{lz}{(1-z)^{l+1}},|z|<1.$
And the renewal equation is given by
\begin{eqnarray}
m_{\alpha,l}(t)=l(D_{0+}^{-\alpha l}\phi_{\alpha})(t)+
(m_{\alpha,l}*\phi_{\alpha})(t)
\end{eqnarray}
We can derive the following from  the expression $\displaystyle
\widetilde{p}_{\alpha,l}(k,t)=t^{\alpha l-1}E_{\alpha,\alpha l}^{l}
(-(1-e^{-k})t^{\alpha})$:
\begin{equation}
 m_{\alpha,l}(t)=-\frac{\partial}{\partial
k}\widetilde{p}_{\alpha,l}(k,t)|_{k=0}
=\frac{lt^{l\alpha+\alpha-1}}{\Gamma((l+1)\alpha)}.
\end{equation}
{\bf\large  Fractional Differential-Difference Equations}\\
And we also get the following analogous differential-difference Equations
\begin{eqnarray}
 \widetilde{P}_{n,\alpha,l}(s)&=&\frac{(l)_{n}}{n!(1+s^{\alpha})^{
l+n}}
\\
\frac{ns^{\alpha}}{n+l-1} \widetilde{P}_{n,\alpha,l}(s)
&=&\widetilde{P}_{n-1,\alpha,l}(s)
-\frac{n}{n+l-1}\widetilde{P}_{n,\alpha,l}(s)
\\
 P_{0,\alpha,l}(t)&=&t^{\alpha
l-1} E_{\alpha,\alpha l}^{l}(-t^{\alpha})
\\
~\frac{n}{n+l-1}D_{0+}^{\alpha}P_{n,\alpha,l}(t)&=&P_{
n-1,\alpha,l}(t)-\frac{n}{n+l-1}P_{n,\alpha,l}(t)
\end{eqnarray}
{\large\bf Erlang Densities}\\
For the corresponding Erlang densities, we have
\begin{eqnarray}
q_{n,\alpha,l}(t)&=&q_{n,\alpha}^{*l}(t)
= t^{ln\alpha-1} E_{\alpha,l\alpha
n}^{ln} (-t^{\alpha})
\\
&=& \frac{t^{ln\alpha-1}}{(ln-1)!}
\sum_{k=0}^{\infty}\frac{\Gamma(ln+k)}{\Gamma( \alpha
(ln+k) )} \frac{(-1)^{k}t^{\alpha k}}{k!}
\end{eqnarray}
{\large\bf Fractional Integral Equation of the Cotinuous Time Random Walk}\\
The fractional integral equation of the CTRW is as
follows:
\begin{eqnarray*}
 \widetilde{\widetilde{p}}_{\alpha,l}(k,s)&=&\left(\frac{
\widetilde{\Psi}_{\alpha}
(s) } {1-\widetilde{\phi}_{\alpha}(s)\widetilde{w}(k)}\right)^{l}
\\
\left(\sum_{j=0}^{l}{l\choose j}(-1)^{j}(
\widetilde{\phi}_{\alpha}(s)\widetilde{w}(k))^{j}\right)\widetilde{\widetilde{p}
}_{\alpha,l}(k, s) 
&=&\left(\widetilde{\Psi}_{\alpha} (s)  \right)^{l}
\\
\sum_{j=0}^{l}{l\choose j}(-1)^{j}
\widetilde{\phi}_{\alpha}^{j}(s)\widetilde{w}^{j}(k)\widetilde{\widetilde{p}}_{
\alpha,l}(k, s)
&=&\left(\widetilde{\Psi}_{\alpha}
(s)  \right)^{l}
\\
p_{\alpha,l}(x,
t)=\Psi_{\alpha}^{*l}
(t)\delta(x)&-&\sum_ {
j=1}^{ l}{
l\choose j}(-1)^{j} \phi^{*j}_{\alpha}(t)*p_{\alpha,l}(x,
t)*w^{*j}(x).
\end{eqnarray*}
{\large\bf Fractional Version of Kolmogorov-Feller Equation}
\begin{eqnarray}
 \widetilde{\widetilde{p}}_{\alpha,l}(k,s)&=&\left(\frac{
\widetilde{\Psi}_{\alpha}
(s) } {1-\widetilde{\phi}_{\alpha}(s)\widetilde{w}(k)}\right)^{l}
\\
 (1+s^{\alpha})\widetilde{\widetilde{p}}_{\alpha,l}(k,s)&=& \frac
{\widetilde{\Psi}_{\alpha}^{l-1}(s) }
{(1-\widetilde{\phi}_{\alpha}(s)\widetilde{w}(k))^{l}}
\\
 s^{\alpha}\widetilde{\widetilde{p}}_{\alpha,l}(k,s)&=&
-\widetilde{\widetilde{p}}_{\alpha,l}(k,s)+\widetilde{\Psi}_{\alpha}(s)
\widetilde{\widetilde{p}}_{\alpha,l}(k,s)
\\
(D_{0+}^{\alpha} p_{\alpha,l}(x,t))(t)&=& 
-p_{\alpha,l}(x,t)
+\int_{0}^{t}\Psi_{\alpha}(t-t_{1})p_{\alpha,l}(x,t_{1})dt_{1}.
\end{eqnarray}
\begin{eqnarray}
 \widetilde{\widetilde{p}}_{\alpha,l}(k,s)&=&\left(\frac{
\widetilde{\phi}_{\alpha}
(s) } {1-\widetilde{\phi}_{\alpha}(s)\widetilde{w}(k)}\right)^{l}
\\
 (1+s^{\alpha})\widetilde{\widetilde{p}}_{\alpha,l}(k,s)&=& \frac{1 }
{(1-\widetilde{\phi}_{\alpha}(s)\widetilde{w}(k))}\widetilde{\widetilde{p_{
\alpha} }}^{l-1}(k,s)
\\
 (1+s^{\alpha})\widetilde{\widetilde{p}}_{\alpha,l}(k,s)&=& \left(\frac{1 }
{(1-\widetilde{\phi}_{\alpha}(s)\widetilde{w}(k))}-1+1\right)
\widetilde{\widetilde{p_{
\alpha} }}^{l-1}(k,s)
\\
 (1+s^{\alpha})\widetilde{\widetilde{p}}_{\alpha,l}(k,s)&=&
\widetilde{w}(k)\widetilde{\widetilde{p_{
\alpha} }}^{l}(k,s)+{\widetilde{p_{
\alpha} }}^{l-1}(k,s)
\\
 s^{\alpha}\widetilde{\widetilde{p}}_{\alpha,l}(k,s)&=&
-\widetilde{\widetilde{p}}_{\alpha,l}(k,s)+\widetilde{w}(k)\widetilde{\widetilde
{p}}_{\alpha,l}(k,s)+{\widetilde{p_{\alpha} }}^{l-1}(k,s)
\\
(D_{0+}^{\alpha} p_{\alpha,l}(x,t))(t)&=& 
-p_{\alpha,l}(x,t)
+\int_{0}^{x}w(x-x_{1})p_{\alpha,l}(x_{1},t)dx_{1} + p_{\alpha}
^{*(l-1)}(x,t).
\end{eqnarray}

\section{Operational Time and Subordination Integrals}
The concept of operational time is introduced in \cite{gm2},
which is an analogous concept of that in operational calculus. Let
\begin{equation}
 p_{1}(y,t_{*})=\sum_{n=0}^{\infty} \frac{t_{*}^{n}}{n!}e^{-t_{*}} \delta(y-n).
\end{equation}
By using the equality 
$\displaystyle
 \int_{0}^{\infty}e^{-at_{*}}dt_{*}=\frac{1}{a} \mbox{ for }a=
s^{\alpha}+1-e^{-k}
$, we can get the following named as the subordination integral, which connects
Poisson process to $\alpha$-fractional Poisson process when the time $t_{*}$
evolves into $t$ with respect to the L\'{e}vy law via Laplace transformation,
\begin{eqnarray}
 \widetilde{\widetilde{p}}_{\alpha}(k,s)&=& \frac{1}{1+s^{\alpha}
-e^{-k}}= \int_{0}^{\infty} e^{-t_{*}(s^{\alpha}+1-e^{-k})}dt_{*}
\\
&=&  \int_{0}^{\infty} e^{-t_{*}(1-e^{-k})}
e^{-t_{*}s^{\alpha}} dt_{*}\label{560}
\\
p_{\alpha}(x,t)&=& \int_{0}^{\infty} p_{1}(x,t_{*})p_{2}(t_{*},t)
dt_{*}
=\int_{0}^{\infty} \sum_{n=0}^{\infty} \frac{t_{*}^{n}}{n!}e^{-t_{*}}
\delta(x-n)  
\frac{\alpha t_{*}}{t^{1+\alpha}}M_{\alpha}\left(\frac{t_{*}}{t^{\alpha}}\right)
dt_{*}.\label{561}
\end{eqnarray}
{\bf\large Evolution equation for the density $p_{2}(t_{*},t)$ of
$t=t(t_{*})$}\\
In \cite{gm3}, the evolution equation for $p_{2}(t_{*},t)$ is defined by
\begin{eqnarray}
 \widetilde{\widetilde{p_{2}}}(s_{*},s)&=&\frac{1}{s^{\alpha}+s_{*}},
~{\widetilde{p_{2}}}(t_{*},s)=\exp(-t_{*}s^{\alpha}),
~{\widetilde{p_{2}}}(s_{*},t)=t^{\alpha-1}E_{\alpha,\alpha}(-s_{*}t^{\alpha})
\\
s_{*}\widetilde{\widetilde{p_{2}}}(s_{*},s)-1&=&
-s^{\alpha}\widetilde{\widetilde{p_{2}}}(s_{*},s)
\\
\frac{\partial}{\partial t_{*}} p_{2}(t_{*},t) &=& -(D_{0+}^{\alpha}
p_{2}(t_{*},t))(t),~~p_{2}(0+,t)=\delta(t)
\end{eqnarray}
In \cite{mpg1}, the famous initial value problem, known as signalling problem,
with the partial differential equation $\displaystyle \frac{\partial}{\partial
t} p(t_{*},t)= \frac{\partial^{2}}{\partial t_{*}^{2}}p(t_{*},t)$ and the
initial conditions, $\displaystyle \lim_{t_{*}\rightarrow 0+}
p(t_{*},t)=\delta(t),t>0,
\lim_{t_{*}\rightarrow +\infty}p(t_{*},t)=0$, has the
following Laplace transform of the solution and the solution
\begin{equation}
 \widetilde{p}(t_{*},s)=e^{-t_{*}s^{\frac{1}{2}}}
\mbox{ and }
p(t_{*},t)=\frac{t_{*}}{2\sqrt{\pi}t^{\frac{3}{2}}}\exp\left(-\frac
{t_{*}^{2}}{4t}\right), 0<t<\infty .
\end{equation}
Therefore this partial differential equation can describe the behavior of the
function $p(t_{*},t)$, which is the integrand of \eqref{561}, or the
relationship between two different time lines
for the case when $\alpha=0.5$.
In \cite{gm2}, the authors have found $q_{\alpha}(t_{*},t)$ for their processes
and gave a relation with $p_{2}(t_{*},t)$.
\begin{eqnarray}
 q_{\alpha}(t_{*},t)&=&t^{-\alpha}M_{\alpha}(t_{*}t^{-\alpha})
\longleftrightarrow \widetilde{q_{\alpha}}(t_{*},s) =s^{\alpha-1}
e^{-t_{*}s^{\alpha}} 
\\
p_{2}(t_{*},t)&=&\frac{\alpha
t_{*}}{t^{1+\alpha}}M_{\alpha}\left(\frac{t_{*}}{t^{\alpha}}\right)
\longleftrightarrow \widetilde{p_{2}}(t_{*},s)=e^{-t_{*}s^{\alpha}}
\\
q_{\alpha}(t_{*},t) &=&(I_{0}^{1-\alpha}p_{2}(t_{*},t_{1}))(t)
\end{eqnarray}
From the above relation, the processes developed here evolve in a different
time line compared with the processes doing in a certain time
fashion in \cite{gm2}. \\
Furthermore, we have similar results for the cases in section 5.2.
\begin{eqnarray}
 \widetilde{\widetilde{p}}_{\alpha,l}(k,s)&=& \underbrace
{\left(\frac {
\widetilde { \phi_ {\alpha} (s)}}{1-\widetilde{\phi_{\alpha}
(s)}\widetilde{w(k)}}\right)\cdots \left(\frac{
\widetilde { \phi_ {\alpha} (s)}}{1-\widetilde{\phi_{\alpha}
(s)}\widetilde{w(k)}}\right) }_{l \mbox{ times}}
\\
&=&
\underbrace
{\left( \frac{1}{1+s^{\alpha}
-e^{-k}} \right)\cdots  \left( \frac{1}{1+s^{\alpha}
-e^{-k}} \right)}_{l \mbox{ times}}
\\
\nonumber &=&
  \underbrace{\int_{0}^{\infty} e^{-t_{*,1}(1-e^{-k})}
e^{-t_{*,1}s^{\alpha}} dt_{*,1}\cdots \int_{0}^{\infty}
e^{-t_{*,l}(1-e^{-k})}
e^{-t_{*,l}s^{\alpha}} dt_{*,l} }_{l
\mbox{ times}}
\end{eqnarray}
\begin{eqnarray}
\nonumber&=&
  \underbrace{\int_{0}^{\infty}\cdots\int_{0}^{\infty}}_{l
\mbox{ times}}\exp(-(1-e^{-k})(t_{*,1}+t_{*,2}+ \cdots +t_{*,l}))
\\
&&\times \exp(-s^{\alpha}(t_{*,1}+t_{*,2}+ \cdots +t_{*,l}))dt_{*,1}\cdots
dt_{*,l} 
\\
\nonumber &=&
  \underbrace{\int_{0}^{\infty}\cdots\int_{0}^{\infty}}_{l
\mbox{ times}}\sum_{n=0}^{\infty} \frac{(t_{*,1}+t_{*,2}+ \cdots
+t_{*,l})^{n}e^{-nk}}{n!} \exp(-(t_{*,1}+t_{*,2}+ \cdots
+t_{*,l}))
\\
&&\times \exp(-s^{\alpha}(t_{*,1}+t_{*,2}+ \cdots +t_{*,l}))dt_{*,1}\cdots
dt_{*,l} \label{563}
\\
&=&
 \underbrace{\int_{0}^{\infty}\cdots\int_{0}^{\infty}}_{l
\mbox{ times}} p_{*}(x,t_{*,1},t_{*,2}, \ldots ,t_{*,l})
p_{2*}(t_{*,1},t_{*,2}, \ldots ,t_{*,l},t)
dt_{*,1}\cdots dt_{*,l} 
\end{eqnarray}
\begin{equation}
 p_{2*}(t_{*,1},t_{*,2}, \ldots ,t_{*,l},t)= \frac{\alpha
(t_{*,1}+t_{*,2}+ \cdots
+t_{*,l})}{t^{1+\alpha}}M_{\alpha}\left(\frac{(t_{*,1}+t_{*,2}+ \cdots
+t_{*,l})}{t^{\alpha}}\right)
\end{equation}

\section{Statiatical Analysis and Some Remarks}
{\bf\large A Random Variable as a Product of Two Independent Random Variables}\\
In \cite{hms2, ma1}, the Mittag-Leffler random variable is well examined. To
support our theory, we will examine the statistical techniques used there in
detail and $p_{2}(t_{*},s)$ in the subordination integral will be reinterpreted.
\begin{thm}\label{l61}
{\bf\large(The Law of Total Expectation)}
 \begin{equation}
  E(x)=E[E(x|y)]
 \end{equation}
where $x$ and $y$ are random variables having a joint distribution, all the
expectations do exist and in the marginal space of $y$, the outside expectation
is calculated.
\end{thm}
In \cite{ma1}, the author statistically describes the structure of
Mittag-Leffler random variable as a product of two
independent random variables by using the lemma \ref{l61}. The method is as
follows:
Let $t_{*}$ be a exponential random variable with the density function
$g_{1}(t_{*})=e^{-t_{*}},~0\leq t_{*}<\infty$ and $y$ be a positive L\'{e}vy
random variable with the density function $\displaystyle g_{2}(y)= \frac{1}{2
\pi i}\oint_{L} \frac{\Gamma(\frac{1}{\alpha} - \frac{s}{\alpha}  )}{\alpha
\Gamma(1-s) }y^{-s}ds, ~~ 1>Re(s)>0,~\alpha>0,y>0$. Then the Mittag-Leffler
random variable $t$ has the form of the product of two independently
distributed random variables $\displaystyle t_{*}^{\frac{1}{\alpha}}, y$.
The proof goes as follows using expectations and their properties below: 
\begin{itemize}
 \item $\displaystyle E[e^{-sy}]=e^{-s^{\alpha}}$
\item $\displaystyle E[\exp(-\{st_{*}^{\frac{1}{\alpha}}\}y)|t_{*}]=
e^{-s^{\alpha}t_{*}}$
\item $\displaystyle E[E[\exp(-\{st_{*}^{\frac{1}{\alpha}}\}y)|t_{*}]]=
E[e^{-s^{\alpha}t_{*}}] =\int_{0}^{\infty} e^{-t_{*}} e^{-s^{\alpha}t_{*}}
dt_{*} =\frac{1}{1+s^{\alpha}} =E[e^{-st}]$.
\end{itemize}
Assume that there exists the random variable $1$ of the constant function
$g(x)=1$ in the ordinary space. As analysed above, the Mittag-Leffler variable
$u$ of $\phi(u)=
u^{\alpha-1}E_{\alpha,\alpha}(-u^{\alpha})$ can be considered as the product of
two independent random variables $t_{1}$ exponentially distributed and $y$
L\'{e}vy-distributed with index $\alpha$, 
namely $\displaystyle u=yt_{1}^{\frac{1}{\alpha}}$. Note that this technique is
well known for real density functions and J-transform will take the place
for non-density functions which is magic and a different
recipe of Mellin convolution technique. We shall analyse the relationship
between a density function $f(t_{1})$
and its cumulative distribution $\Phi(t_{2})$ statistically. Look at the
mathematical
relation between them in the ordinary space  in the below form
\begin{equation}
 \Phi(t_{2})=\int_{0}^{t_{2}}1 \cdot
f(t_{1})dt_{1}=\int_{0}^{t_{2}}g(t_{2}-t_{1})
\cdot f(t_{1})dt_{1}.\label{52}
\end{equation}
This relation can be interpreted in statistics as follows.
We can consider $\Phi(t_{2})$ as the density function of the random variable
$t_{2}=t_{1}+1$ which is the Laplace convolution of two functions as in
(\ref{52}).
With the help of J-transformation, we have the corresponding generalized
random
variable $y_{1}=y 1^{\frac{1}{\alpha}}$ of $1$ in the $\alpha$-fractional space
with the generalized density function $\displaystyle \frac{y_{1}^{\alpha
-1}}{\Gamma(\alpha)}$. Because there are no statistical methods for this
correspondence except J-transformation. From the above analysis, we
take the random variable $\displaystyle t_{3}=yt_{1}^{\frac{1}{\alpha}} +y
1^{\frac{1}{\alpha}}$ for the marginal cumulative distribution of the density
function $\phi_{\alpha}(u)$ of the Mittag-Leffler random variable
$u=yt_{1}^{\frac{1}{\alpha}}$.
That leads to the following relation between a density function and its
marginal cumulative distribution in the $\alpha$-fractional space
\begin{equation}
\Phi_{\alpha}(t_{3})=\int_{0}^{t_{3}}\frac{(t_{3}-u)^{\alpha-1}}{\Gamma(\alpha)}
 u^{\alpha-1}E_{\alpha,\alpha}(-u^{\alpha})du 
\end{equation}
which supports the theory developed in section \ref{5}.\\
{\bf\large The Subordination Integral}
\begin{cor}\label{c71}
Let $h_{1}(x)$ be a density function of a random variable $x$. Let $h_{2}(y|x)$
be a conditional density function of a variable $y$ given $x$. Define
\begin{equation}
 h(x,y)=\begin{cases}
         h_{2}(y|x)h_{1}(x) \mbox{ if } h_{1}(x)>0 \mbox{ for all }x \\
	0 \mbox{ if } h_{1}(x)=0.
        \end{cases}
\end{equation}
Then $h(x,y)$ is a joint density function and $h_{1}(x)$ is the marginal
density function of $x$.
\end{cor}
We have the subordination integral by
\begin{eqnarray}
 \widetilde{\widetilde{p}}_{\alpha}(k,s)&=&   \int_{0}^{\infty}
e^{-t_{*}(1-e^{-k})} e^{-t_{*}s^{\alpha}}
dt_{*}, ~p_{\alpha}(x,t) =
\int_{0}^{\infty} p_{1}(x,t_{*})p_{2}(t_{*},t)dt_{*}.
\end{eqnarray}
From the Corrollory \ref{c71}, $p_{3}(t_{*},t)=e^{-t_{*}}p_{2}(t_{*},t)$
is  a joint distribution of $t_{*}$ and $t$. Therefore we get, via statistical
paths, 
 \begin{eqnarray}
 \int_{0}^{\infty} e^{-t_{*}} e^{-s^{\alpha}t_{*}} dt_{*}
&=&\frac{1}{1+s^{\alpha}} =\widetilde{\phi_{\alpha}}(s) \Longrightarrow
\phi(t):\mbox{ the density in the marginal space of }t
\\
 \int_{0}^{\infty} e^{-t_{*}} p_{2}(t_{*},t) dt&=&e^{-t_{*}}:\mbox{ the density
in the marginal space of }t_{*}.
 \end{eqnarray}
The evolution equation for $p_{3}(t_{*},t)$ can be derived to
\begin{equation}
p_{3}(t_{*},t)=-\frac{\partial}{\partial t_{*}} p_{3}(t_{*},t) 
-(D_{0+}^{\alpha}p_{3}(t_{*},t))(t),~p_{3}(0+,t)=\delta(t).
\end{equation}
{\bf\large Stochastic Analysis and Pathways}\\
Let $f(t_{*})=e^{-t_{*}}$ and $\phi_{\alpha}(t)$ be density functions of
$t_{*}$ and $t$. Then we get the following stochastic processes:
\begin{eqnarray}
 f(t_{*})=e^{-t_{*}}&\Longrightarrow& \lim_{l\rightarrow \infty}
lf^{*l}(lt_{*})=\delta(t_{*}-1),~ \mathcal{L}\{\delta(t_{*}-1)\}(s)=e^{-s}
\\
\phi_{\alpha}(t)=t^{\alpha-1}E_{\alpha,\alpha}(-t^{\alpha})
&\Longrightarrow& \lim_{l\rightarrow \infty} l^{\frac{1}{\alpha}}
\phi_{\alpha}^{*l}(l^{\frac{1}{\alpha}}t)=\delta_{\alpha}(t -1),~
\mathcal{L}\{\delta_{\alpha}(t-1)\}(s)=e^{-s^{\alpha}}\label{789}
\\
F(t)=f^{*n}(t_{*}) &\Longrightarrow& \lim_{l\rightarrow \infty}
lF^{*l}(lt_{*})=\delta(t_{*}-n),~ \mathcal{L}\{\delta(t_{*}-n)\}(s)=e^{-ns}
\\
\psi_{\alpha}(t)=\phi_{\alpha}^{*n}(t) &\Longrightarrow& \lim_{l\rightarrow
\infty} l^{\frac{1}{\alpha}}
\psi_{\alpha}^{*l}(l^{\frac{1}{\alpha}}t)=\delta_{\alpha}(t
-n^{\frac{1}{\alpha}}),~
\mathcal{L}\{\delta_{\alpha}(t-n^{\frac{1}{\alpha}})\}(s)=e^{-ns^{\alpha}}
\end{eqnarray}
\begin{eqnarray}
  f(t_{*})=e^{-t_{*}}&\Longrightarrow& \lim_{l\rightarrow
\infty} \frac{l}{n}f^{*l}(\frac{l}{n}t_{*})=\delta(t_{*}-n),~
\mathcal{L}\{\delta(t_{*}-n)\}(s)=e^{-ns}\label{790}
\\
\phi_{\alpha}(t)=t^{\alpha-1}E_{\alpha,\alpha}
(-t^{\alpha } )
&\Longrightarrow& \lim_{l\rightarrow \infty}
\left(\frac{l}{t_{*}}\right)^{\frac{1}{\alpha}}
\phi_{\alpha}^{*l}(\left(\frac{l}{t_{*}}\right)^{\frac{1}{\alpha}}t)=\delta_{
\alpha } (t-(t_{*})^{\frac{1}{\alpha}}),\label{795}
\\
&&\mathcal{L}\{\delta_{\alpha}(t-(t_{*})^{\frac{1}{\alpha}})\}(s)=e^{-t_{*}s^{
\alpha}}\label{796}.
\end{eqnarray}
The positive-oriented (extreme) stable process $t=t(t_{*})$ with
$p_{2}(t_{*},t)$ is described with $t_{*}=n\tau_{*}$ and $\tau_{*}$ as a
step-size, $\tau_{*}>0$, in \cite{gm3}, which is \eqref{796} with
$t_{*}=n\tau_{*}, ~n=1,2,3,\ldots$. \eqref{789} is treated as a stochastic
process by applying
Laplace technique without the explcit form in \cite{pi1} and appears explicitly
in
\cite{h1}. The explcit form in \cite{h1} appears firstly in the literature. But
we emphasise that in \cite{pi1}, the process was not treated in
the $\alpha$-fractional space.\\ 
In the last remark of \cite{pi1}, the author describes the Mittag-Leffler
process in \cite{pi1} as a stochastic process subordinated to a stable process
by directing the gamma process, or randomizing the parameter $s$ with gamma
distribution. For the each $l$-th event of the time line of the Poisson process
in this paper,
\begin{equation}
 p_{4}(t_{*},t)=f^{*l}(t_{*})p_{2}(t_{*},t)
\end{equation}
\begin{eqnarray}
 \int_{0}^{\infty} f^{*l}(t_{*}) e^{-s^{\alpha}t_{*}} dt_{*}
&=&\left(\frac{1}{1+s^{\alpha}}\right)^{l} =\widetilde{\phi_{\alpha}}^{l}(s)
\Longrightarrow
\phi_{\alpha}^{*l}(t):\mbox{ the density in the marginal space of }t \nonumber
\\
 \int_{0}^{\infty} f^{*l}(t_{*}) p_{2}(t_{*},t) dt&=&f^{*l}(t_{*}):\mbox{ the
density
in the marginal space of }t_{*}.\nonumber
 \end{eqnarray}
{\bf\large J-transformation Correspondence}\\
The following table from the author's paper \cite{h1} supplies correspondence
between ordinary space and $\alpha$-fractional space via J-
transformation:\\
\textbf{\tiny\bf J-transformation CORRESPONDENCE}\\
\begin{tabular}{|c|c|c|}
  \hline
  \mbox{ordinary space} &  & \mbox{$\alpha$-level(fractional) space}  \\
  \hline
  $1$ &  & $\frac{t^{\alpha-1}}{\Gamma(\alpha)}$ \\
\hline $t$&&$\frac{t^{2\alpha-1}}{\Gamma(2\alpha)}$\\
  \hline
  $ e^{-t}$ &  & $t^{\alpha-1}E_{(\alpha,\alpha)} (- t^{\alpha}) $ \\
  \hline
$\frac{t^{l-1}}{(l-1)!}e^{-t}$ && $t^{\alpha l-1}E_{\alpha,\alpha l}
^{l}(-t^{\alpha})$ \\
\hline
\end{tabular}\\
In this paper, we develop the theory of $\alpha$-fractional Poisson process. As
stated in section 3, the renewal function of the
Poisson process is $m(t)=t$, but that of the $\alpha$-fractional Poisson process
is $\displaystyle m(t)=\frac{t^{2\alpha-1}}{\Gamma(2\alpha)}, ~~0<\alpha<1$. The
renewal function in the $\alpha$-fractional world looks peculiar. On the other
hand, we get the same function corresponding from the renewal function of
ordinary Poisson case in the J-transformation correspondence table above, namely
$\displaystyle t \longleftrightarrow \frac{t^{2\alpha-1}}{\Gamma(2\alpha)}$. 
It is interesting that the J-transformation correspondence table provides the
densities
of two processes as well as their renewal functions.

\section{Remark}
We think that our theory is related to the theory of infinite field extention
and a view of irrationality in mathematics and the theory of 
long-term-memory type in statistics. We also recognise the fact, which should be
overcome by some logical and practical argument, that the equation (\ref{k100})
gives some wierd difference from the usual statistical approach. 
\section{Apendix}
Pochammer symbol is
defined as
\begin{equation}\label{} 
(b)_{k}=b(b+1)\cdots (b+k-1), ~(b)_{0}=1.
\end{equation}

\begin{de}
Let $f(x) \in L(a,b), \alpha \in \mathbb{C}, Re(\alpha)>0$, then
\begin{equation}
(I_{0+}^{\alpha}f)(x)=(D_{0+}^{-\alpha}f)(x)= (~_{0}I_{x}^{\alpha}f)(x)=
(~_{0}D_{x}^{-\alpha}f)(x) =\frac{1}{\Gamma(\alpha)}
\int_{0}^{x} (x-t)^{\alpha-1}f(t)dt, x>0,
\end{equation}
 which is called the Riemann-Liouville left-sided fractional integral of
order $\alpha$ in \cite{mhs1}  or Abel integral operator for
$~_{0}D_{x}^{-\alpha}$ in \cite{gv1}.
\end{de}
\begin{de}
 The Weyl integral of order $\alpha$ is defined as follows:
\begin{equation}
 (~_{x}W_{\infty}^{\alpha}f)(x)=\frac{1}{\Gamma(\alpha)} \int_{x}^{\infty}
(t-x)^{\alpha-1}f(t)dt, (0<x< \infty)
\end{equation}
where $\alpha \in \mathbb{C}, Re(\alpha)>0.$
\end{de}

\begin{de}
 The Mainardi function is defined as
\begin{equation}
M_{\alpha}(z):=\sum_{n=0}^{\infty} \frac{(-z)^{n}}{n!\Gamma[-\alpha
n+(1-\alpha)]}= \frac{1}{\pi}\sum_{n=1}^{\infty}\frac{(-z)^{n-1}}{(n-1)!}
\Gamma(\alpha n)sin(\pi \alpha n) 
\end{equation}
with $z \in \mathbb{C}$ and $0<\alpha<1$. For detailed information, we refer to
\cite{m1}.
\end{de}

\begin{de}
 Mittag-Leffler function of 3 parameters is defined by
\begin{equation}
 E_{\alpha,\beta}^{\gamma}(x)=\sum_{k=0}^{\infty} \frac{(\gamma)_{k} x^{k}}{
\Gamma(\alpha k + \beta) k!} (\alpha, \beta, \gamma\in \mathbb{C};Re(\alpha)>0,
Re(\beta)>0).
\end{equation}
\end{de}
The importance and usage of the Mittag-Leffler function have been enlarged
immensely, especially in the area of Fractional Calculus. For theoretical
approachs and its applications,  there are many good books and papers dealing
with functions of Mittag-Leffler type as can be found in the reference.

\vskip 0.2cm
\textbf{\emph{Laplace Transformation Correspondence}}\\
\begin{tabular}{|c|c|c|}
  \hline
function &$\longleftrightarrow$  & Laplace transform \\
 \hline
 $\displaystyle E_{\alpha}(at^{\alpha})=\sum_{n=0}^{\infty} \frac{(-1)^{n}
(at^{\alpha})^{n}}{\Gamma(\alpha n+1)}$&&$ \displaystyle\frac{\frac{s^{\alpha
-1}}{a}}{1+\frac{s^{\alpha}}{a}}$
\\
 \hline
$\displaystyle\frac{t^{n\alpha}}{n!}E^{(n)}_{\alpha}(-t^{\alpha}
)$&&$\displaystyle \frac{s^{ \alpha-1}}{ (1+s^{\alpha})^{n+1}} $
\\
 \hline
$\displaystyle
\frac{x^{-n\alpha-1}}{\Gamma(-n\alpha)}$&&$\displaystyle s^{n\alpha}$\\
 \hline
$\displaystyle
\frac{\alpha n}{t^{\alpha+1}}M_{\alpha}
\left(\frac{n}{t^{\alpha}} \right)$ in \cite{gm1,m1}
&&$\displaystyle
e^{-ns^{\alpha}}$ 
\\
\hline
$\displaystyle \lim_{\gamma\rightarrow\infty}
t^{\alpha\gamma-1}\left(\frac{\gamma}{n}\right)^{\gamma}E^{\gamma}_{(\alpha,
\alpha\gamma) }
\left(-\frac{\gamma t^{\alpha}}{n}\right)$ in \cite{h1}
&&$\displaystyle
e^{-ns^{\alpha}}$ 
\\
 \hline
$\displaystyle t^{\alpha k+\beta-1}  
E^{(k)}_{(\alpha,\beta)}(-at^{\alpha}),E_{\alpha,\beta}^{(k)}(y)=\frac{d^{k}}{dy
^{k}}E_{(\alpha,\beta)}(y) $&&$\displaystyle\frac { k!s^ {
\alpha-\beta } } { \left(a+s^{\alpha}\right)^{-(k+1)}},
~Re(s)>1$
\\
 \hline
$\displaystyle t^{\beta-1}  
E_{(\alpha,\beta)}(-ct^{\alpha})$&&$\displaystyle\frac{s^{\alpha-\beta}}{c+s^{
\alpha }}$
\\
 \hline
  \end{tabular}


\begin{thebibliography}{7}


\bibitem{bc1}
Beck C. and Cohen E.G.D., Superstatistics, Physica A, 322:
267-275, 2003.

\bibitem{c1}
Cox D.R., Renewal Theory, 2-nd Edition, Methuen, London, 1967.

\bibitem{f1}
  Feller, W.  An Introduction to Probability Theory and Its Applications, Vol.
II, Wiley, New York (1966).



\bibitem{gs1}
Gel'fand I.M. and Shilov G.E., Generalized functions, Vol. 1. Academic Press,
New York, 1964.

\bibitem{gm1}
Gorenflo R. and Mainardi F., On the fractional Poisson process and
the discretized stable subordinator, preprint.

\bibitem{gm2}
Gorenflo R. and Mainardi F., Laplace-Laplace analysis of the fractional Poisson
process, arXiv:1305.5473v1 [math.PR] 23 May 2013.

\bibitem{gm3}
Gorenflo R. and Mainardi F., Parametric Subordination in Fractional Diffusion
Processes, in J. Klafter, S.C. Lim and R. Metzler (Editors), Fractional
Dynamics, World Scientific, Singapore, 2012, chapter 10, pp. 229-263. 


\bibitem{gv1}
Gorenflo R. and Vessella S., Abel integral equations: Analysis and
applications, Lecture Notes in Mathematics No. 1461, Springer-Verlag, 1991.


\bibitem{hms2}
Haubold H. J., Mathai A. M. and Saxena R.K., Mittag-Leffler functions and their
applications, Journal of Applied Mathematics (ID 298628, 2011, 51 pages)(2011).

\bibitem{h1}
 Jung Hun  Han, On the Levy density function, arXiv:1102.2709v1 [math.ST] 14 Feb
2011.

\bibitem{h2}
Jung Hun Han, One-sided L\'{e}vy stable distributions, Intellectual Archive
Vol. 1, No. 3, 174-184, July 2012, (arXiv:1102.2713)

\bibitem{h3}
Jung Hun Han, Gamma function to Beck-Cohen superstatistics, Physica
A: Statiatical Mechanics and its applications 392 (2013), pp. 4288-4298.

\bibitem{kst1}
Kilbas A.A., Srivastava H.M. and Trujillo J.J., Theory and Applications of
Fractional Differential Equations, Elsevier, Amsterdam, 2006. 

\bibitem{m1}
Mainardi F., Fractional Calculus and Waves in Linear Viscoelasticity, Imperial
College Press, London, 2010.

\bibitem{mgs1}
 Mainardi F., Gorenflo R. and Scalas, A fractional generalization of the
Poisson processes, Vietnam J. of Math., 32, SI (2004), 53-64.

\bibitem{mlg1}
 Mainardi F.,  Luchiko Y. and Gorenflo R. , The fundamental solution of the
space-time fractional diffusion equation,Fractional calculus and Applied
Analysis, Vol.4,No.2(2001), 153-192.


\bibitem{mpg1}
Mainardi F., Paradisi P and Gorenflo R., Probability distributions
generated by fractional diffusion equations, FRACALMO preprint,
www.fracalmo.org, 1997. To appear in: J. Kertesz and I. Kondor (editors):
Econophysics. Kluwer Academic Publishers, Dordrecht (NL), 1998, pp 39ff.
ALTHOUGH PLANNED, the book did NEVER appear!

\bibitem{ma1}
 Mathai A.M., Some properties of Mittag-Leffler functions and matrix-variate
analogues: a statistical perspective, Fractional Calculus $\&$ Applied Analysis,
3(1), 113-132 (2010).




\bibitem{mh1}
Mathai A.M. and  Haubold Hans J., Mittag-Leffler functions to
Pathway Model to Tsallis statistics, Integral Transforms and Special Functions,
21(11), 867-875(2011).

\bibitem{mh2}
Mathai A.M. and  Haubold Hans J., Special Functions for Applied
Scientists, Springer, New York, 2008.


\bibitem{mhs1}
Mathai A.M., Saxena R.K. and  Haubold Hans J., The H-function:
Theory and Applications, Springer, New York, 2010.


\bibitem{pi1}
Pillai R.N., On Mittag-Leffler functions and related distributions, Ann. Inst.
Statist. Math., 42(1990), 157-161.

\bibitem{p1}
Podlubny I., Fractional Differential Equations, Academic Press, San Diego, 1999.



\bibitem{t1}
Tsallis, C., Introduction to Nonextensive Statistical Mechanics:
Approaching a complex world, Springer, New York, 2009.


\bibitem{w1}
Weiss G.H., Aspects and applications of random walks, North-Holland, Amsterdam,
1994.



\end{thebibliography}
\end{document}